\documentclass[a4paper,10pt]{article}
\usepackage[utf8]{inputenc}
\usepackage{amsmath,epsfig}
\usepackage{amssymb,amsfonts}
\DeclareMathOperator{\Emb}{\mathrm{Emb}}
\DeclareMathOperator{\Coh}{\mathrm{Coh}}

\title{The Critical Numbers of Rankin-Selberg Convolutions of Cohomological Representations}
\author{Claus G\"unther Schmidt}

\begin{document}

\maketitle

\begin{abstract}
We study the critical numbers of the Rankin-Selberg convolution of arbitrary pairs of cohomological cuspidal automorphic representations and  we parametrize these critical numbers by  certain 1-dimensional subrepresentations 
attached to the corresponding pair of finite dimensional representations of the related general linear groups.

\end{abstract}

\section{Introduction}

For arbitrary natural numbers $n$ and $m$ let $\pi$ and $\sigma$ denote cuspidal automorphic representations of $GL_n(\Bbb A)$ and $GL_m(\Bbb A)$ respectively over the adele
ring $\Bbb A$ of $\Bbb Q$. Jacquet, Piatetski-Shapiro and Shalika [4] introduced for such pairs $(\pi ,\sigma )$ an $L$-function $L(\pi ,\sigma ,s)$ which up to a few special cases
is an entire function of the variable $s$. In analogy with Deligne's notion of {\it critical values} of motivic $L$-functions [2] we would like to study the values of 
$L(\pi ,\sigma ,s)$ at {\it critical numbers} $t \in \frac{n-m}{2} + \Bbb Z$ and in particular their arithmetic properties. Assuming that $\pi$ and $\sigma$ are cohomological in 
the case $m=n-1$ these critical values are quite well understood by [5]. See also Januszewski's contributions [6] for totally real number fields. The assumption says that there
are finite-dimensional irreducible rational representations $M_\mu$ and $M_\nu$ of $GL_n$ and $GL_m$ respectively of heighest weights $\mu$ and $\nu$ with a certain purity property
such that for the infinity components $\pi _\infty$ and $\sigma _\infty$ the representations $\pi _\infty \otimes M_\mu$ and $\sigma _\infty \otimes M_\nu$ have non-trivial relative Lie algebra
cohomology, i.e. we have
$$H^\bullet (gl_n,K_{n,\infty};\pi_\infty \otimes M_{\mu ,\Bbb C}) \neq 0$$
and
$$H^\bullet (gl_m,K_{m,\infty};\sigma_\infty \otimes M_{\nu ,\Bbb C}) \neq 0,$$
where for a natural number $n$ as usual $gl_n$ denotes the Lie algebra of $GL_n(\Bbb R)$, $K_{n,\infty}=SO_n(\Bbb R) Z_n^+(\Bbb R)$ and $Z_n^+(\Bbb R)$ is the subgroup of matrices of positive determinant in the 
center $Z_n(\Bbb R)$ of $GL_n(\Bbb R)$.
For a given weight $\mu$ the set of representations $\pi$ with this property is usually denoted by $\Coh(\mu)$. The treatment of the critical values for $m=n-1$ relied on the bijection $t \mapsto t + \frac{1}{2}$ in this case between the set $\mathrm{Crit}(\pi _\infty ,\sigma _\infty )$ of
critical numbers and the parameter set $\Emb(\nu ,\check{\mu})$ of  integers $s$ allowing to embed the twists $M_{\nu -s} = \det ^s \otimes M_\nu$ for $\nu -s = (\nu _1 -s,...,\nu _m -s)$ into the 
contragredient $M_{\check{\mu}}$ of $M_\mu$ considered as a $GL_m$-module under the restricted action of $GL_m \hookrightarrow GL_n$ where we had to assume the existence
of at least one such embeddable twist. Eventually this enabled us to interprete the critical values of $L(\pi ,\sigma ,s)$ as resulting from a pairing of appropriate cohomology spaces
with coefficients in $M_\mu$ and $M_\nu$ thus supplying a powerful technique towards algebraicity and integrality statements for the critical values.

  Unfortunately this approach seems not to work in general. Even in the case $m=n-1$ it may happen that $\Emb(\nu ,\check{\mu})$ is empty although critical numbers do exist. In this 
article I will therefore by-pass this obstacle by not only restricting the $GL_n$-action on $M_\mu$ to a suitable smaller subgroup  but also restricting the $GL_m$-action on $M_\nu$
and carefully analyze the correlation between certain irreducible components $M_{\theta _i'(\tilde{\lambda})},\ M_{\theta _i(\tilde{\mu})}$ for $i=1,2$ of modifications of such 
restrictions to certain subgroups $GL_r$ and $GL_{r+1}$. Like in [7] a technical  tool is the ramification law for the restriction of an irreducible $GL_{r+1}$-representation $M_\alpha$ of highest weight
$\alpha$ to $GL_r$ saying that this restriction decomposes as a direct sum of irreducible $GL_r$-representations $M_\beta$ of highest weight $\beta$ each of multiplicity one where $\beta$ varies subject to
the condition 
$$\alpha _\rho \geq \beta _\rho \geq \alpha _{\rho +1}  \mbox{\ \ for \ \ } \rho = 1,...,r.$$
We will  write $\beta \prec \alpha $ in this case and we denote in general by 
$$\Emb (\beta ,\alpha ) := \{s \in \Bbb Z; \beta - s \prec \alpha \}$$
the parameter set of embeddable twists.
The involved modifications in particular incorporate the relative position of the respective Langlands parameters  $(w,l)$ of $\pi _\infty$
and $(w',l')$ of $\sigma _\infty$ via a certain {\it position tupel} $a \in  \Bbb Z^m$ assuming without loss of generality that $l_1>l_1'$. The components of $a$ are uniquely determined by the 
requirement
$$l_{a_j} > l'_j \geq  l_{1+a_j} \mbox{\ \ for \ \ } j=1,...,m$$
and in a first step the dominant weight $\nu$ is replaced by the modified weight $\lambda$ where we put $\lambda _j := \nu _j + a_j -j$.  In the {\it non-exceptional} case where $n$ and $m$ are not both odd   $\lambda $ turns out to be 
a dominant weight with purity property. In the exceptional case a slight further modification leads to  these qualities.

To begin with in Proposition 2.1 we characterize the critical numbers by the inequality 
$$|t - \kappa| < L + \frac{1}{2} $$
for $\kappa := \frac{1}{2} (w + w' + 1)$ and 
$$L:= \frac{1}{2} \min \{|l_i - l'_j|;\mbox{ all } i,j \mbox{ such that } i\neq \frac{n+1}{2} \mbox{ or } j\neq \frac{m+1}{2}\}$$ 
with an additional parity condition in the {\it exceptional} case when $n$ and $m$ both are odd. Then we discuss the so-called {\it jump indices} $j$ where the position sequence $a$ increases in
the sense that $a_{j+1}>a_j$. A simple transformation rule  relating the Langlands parameters $(w,l)$ and $(w',l')$ of $\pi_\infty$ and $\sigma_\infty$ with the corresponding weights $\mu$ and $\nu$
allows us to translate the characterizing inequality for critical numbers $t$ into a system of inequalities of the form
$$\check{\mu }_{a_j} \geq \lambda _j - s \geq \check{\mu }_{1+a_j} \mbox{ for } j=1,...,m \leqno(I_j)$$
where $s=t+\frac{n-m}{2}-1$ with slight modifications in the exceptional case (see Proposition 3.2). 
In terms of the jump indices $j_1,...,j_k$ this inequality system can be expressed in the form
$$\check{\mu}_{a_{j_\kappa}} \geq \lambda _{1+ j_{\kappa -1}} -s \geq ...\geq \lambda _{j_\kappa} -s \geq \check{\mu}_{1+{a_{j_\kappa}}} \leqno(S_\kappa)$$
for $\kappa = 1,...,k+1$ with $j_0=0$ and $j_{k+1}=m$. In fact by the dominance property of $\lambda$ the inequality sequence $(S_\kappa)$ is equivalent to demanding just the first and the last inequality.
This observation is vital for us since it allows by a careful book-keeping to translate the inequality system via the ramification law into an embedding statement of $GL_r$-representations for $r:=k+1$.
This we work out first for a general system of inequalities in Proposition 4.1 which we then use as a guideline for the composition of suitably modified highest  weights $\theta _i'(\tilde{\lambda})$
for $GL_r$ and $\theta _i(\tilde{\mu})$ for $GL_{r+1}$ with purity property for $i=1,2$. Eventually our 
main result in Theorem A (in the non-exceptional case) and Theorem B (in the exceptional case) unconditionally parametrizes the critical numbers as simultaneous parameters  of embeddings 
that belong to $\Emb (\theta _i'(\tilde{\lambda}),\theta _i(\check{\tilde{\mu}}))$ for $i=1,2$ thus
generalizing the specific approach in the case $m=n-1$. In contrast to that special case in this paper the representation $\pi $ a priori does not play a privileged role compared to $\sigma$. 
This fact is taken into account by the more symmetric formulation of the results in the two corollaries attached to the theorems where the critical numbers essentially get parametrized by the common $1$-dimensional pieces  of the
canonically associated $GL_r$-representations 
$$(M_{\theta _i(\tilde{\mu})} \otimes M_{\theta _i'(\tilde{\lambda}}))^{SL_r} \mbox{ \ \ for \ \ } i=1,2$$
with pieces given by the  {\it Tate modules} $T_r(s) := M_{(s,...,s)} = \det ^s$ . The very special case $n = m = 1$ is excluded here since the critical
numbers of Dirichlet-$L$-functions are well understood.

One might interprete the results in this paper as a vague hint towards a possibly extended cohomological construction of {\it modular symbols} attached to the critical values of the $L$-function
$L(\pi ,\sigma ,s)$ for cohomological cuspidal representations $\pi$ and $\sigma$ in general.

\section{Langlands parameters of critical pairs}

As is well known by (3.6) in [8] following Lemme 3.14 in [1] the infinity component $\pi_\infty$ of a cohomological representation $\pi$ of $GL_n$  is an induced representation 
of Langlands type of the form
$$\pi_\infty \cong  J(-w,l) \otimes  sgn^\delta $$
(see also (1.8) in [7]), where the parameters $w\in \Bbb Z$ and $l=(l_i)\in \Bbb Z^n$ belong to the set
$$L^+_0(n) = \{(w,l) \in \Bbb Z \times \Bbb Z^n ; l_i > l_{i+1}, l_i + l_{n+1-i} = 0, w+l_i \equiv n+1 \bmod 2\}$$
and $\delta \in \{0,1\}$. Similarly for a cohomological representation $\sigma$ of $GL_m$ the infinity component $\sigma _\infty$ is of the form
$$\sigma_\infty \cong J(-w',l') \otimes sgn^{\delta '}$$
with $(w',l')\in L_0^+(m)$, where we use the same sign convention for $w$ and $w'$ as in [7](1.8).

In analogy with Deligne's notion of critical values of motivic L-functions we call a number 
$t \in  \frac{m+n}{2} +\Bbb Z$ a {\it critical number} if the archimedean L-function $L(\pi_\infty , \sigma_\infty ,s)$ and its counterpart $L(\check{\pi}_\infty ,\check{\sigma}_\infty ,1-s)$
in the functional equation with the contragredients $\check{\pi}_\infty$ and $\check{\sigma}_\infty$ do not have a pole in $t$.
Recall that  the L-function $L(\pi_\infty ,\sigma_\infty ,s)$ is defined to be the L-function of the representation $\tau := \pi_\infty^W \otimes \sigma_\infty ^W$ of the  Weil group $W_{\Bbb R}$
where $\pi_\infty ^W$ and $\sigma_\infty ^W$ denote the semisimple representations of $W_{\Bbb R}$ attached to $\pi_\infty$ and $\sigma_\infty$ by the archimedean local Langlands correspondence.
We denote by  $\mathrm{Crit} = \mathrm{Crit}(\pi_\infty ,\sigma_\infty )$ 
the set of critical numbers $t$ and we call the pair $(\pi _\infty ,\sigma _\infty )$ a {\it critical pair}, if $\mathrm{Crit}$  is not empty. This finite set can be described in terms of the associated Langlands parameters
as follows. We put 
$$L_0 := \min \{|l_i - l_j'|; \mbox{ \ all pairs \ } (i,j) \mbox{ such that } i\neq  \frac{n+1}{2} \mbox{ or } j\neq \frac{m+1}{2}\}$$
we set $L:= L_0/2$. Note that the restriction for the choice of $i$ and $j$ under consideration in the minimum 
only matters in the exceptional case $n\equiv m\equiv 1 \mod 2$. Otherwise $i$ and $j$ vary unconditionally. Let $\kappa := \frac{1}{2}(w+w'+1)$ and $\kappa ':= \kappa -\frac{1}{2}$.

\bigskip

{\bf Proposition 2.1:} {\it A number $t \in \frac{m+n}{2} + \Bbb Z$ is critical if and only if 
$$ |t - \kappa | < L + \frac{1}{2} \leqno(I)$$
and in addition in case $n\equiv m\equiv 1 \mod 2$ the number $t$ satisfies the parity condition 
$$ t - \kappa ' \in (2\Bbb N - \epsilon )\cup -(2\Bbb N -1-\epsilon ) \leqno(PC)_\epsilon$$
for $\epsilon = 0,1 $ with $\epsilon \equiv \delta + \delta ' \mod 2$.}

\bigskip

{\bf Proof.} In the special case $m = n-1$ in [7] a proof has been worked out which we now want to extend to the general case. The semisimple representation $\pi_\infty ^W$ of the Weil group $W_{\Bbb R}$
attached to $\pi_\infty$ is of the form
$$\pi_\infty ^W = (l_1,-w/2)\oplus ...\oplus (l_{n/2},-w/2) $$ for even $n$
and
$$\pi_\infty ^W = (l_1,-w/2)\oplus ...\oplus (l_{\frac{n-1}{2}},-w/2)\oplus (sgn^\delta,-w/2) $$ for odd $n$,
where $(l,t)$ and $(\pm ,t)$ for integers $l\geq 1$ and $t\in \Bbb C$ are the irreducible 2-dimensional resp. 1-dimensional representations of $W_{\Bbb R}$ in Knapp's notation ([6], (3.2),(3.3)). 
Analoguously we decompose $\sigma_\infty ^W$ in terms of the parameters $l'_j$, $w'$ and $\delta '$. The L-function of $\tau =\pi_\infty ^W \otimes \sigma_\infty ^W$ is given by the product of the L-functions 
$L(s,\varphi )$ of the irreducible constituents $\varphi$ of $\tau$. Using Deligne's terminology 
$\Gamma_{\Bbb R}(s) := \pi ^{-s/2} \Gamma (s/2)$ and $\Gamma_{\Bbb C}(s) := 2 (2\pi )^{-s} \Gamma (s)$ we can express the $L(s,\varphi )$ by (3.6) in [6] as
$$L(s,(l,t)) = \Gamma_{\Bbb C}(s + t + \frac{l}{2}) \leqno(2.1)$$
and
$$L(s,(sgn^\epsilon ,t)) = \Gamma_{\Bbb R}(s + t + \epsilon ) \leqno(2.2)$$
with $\epsilon = 0,1$.

In order to write down the L-function for $\tau$ (and the contragredient $\check{\tau}$) we must work out the decomposition into irreducible parts. To simplify notation we will sometimes
use the reducible representations $(0,t) := (+,t)\oplus (-,t)$. There are three possible parity combinations for $n$ and $m$ to be treated separately.
\bigskip

{\bf Lemma 2.1:} {\it For $n\equiv m\equiv 0 \mod 2$ we have
$$\tau = \bigoplus _{i=1}^{n/2} \bigoplus _{j=1} ^{m/2} (l_i+l'_j,-\kappa ')\oplus (|l_i-l'_j|,-\kappa ')$$
hence
$$L(s,\tau ) = \prod _{i,j} \Gamma_{\Bbb C}(s-\kappa '+ \frac{l_i+l'_j}{2})$$
$$\cdot \prod _{i,j; l_i\neq l'_j} \Gamma _{\Bbb C}(s-\kappa '+\frac{|l_i-l'_j|}{2})\cdot
\prod _{i,j; l_i=l'_j} \Gamma_{\Bbb R}(s-\kappa ') \Gamma_{\Bbb R}(s-\kappa '+1).$$}
\bigskip 

The proof of Lemma 2.1 is an easy exercice using the rules for the tensor product of the envolved irreducible representations 
$$ (sgn^\epsilon ,t)\otimes (sgn^{\epsilon '},t') = (sgn^{\epsilon +\epsilon '},t+t'),$$
$$ (l,t)\otimes (sgn^\epsilon ,t) = (l,t+t'),$$
$$ (l,t)\otimes (l',t') = (l+l',t+t')\oplus (|l-l'|,t+t')$$
as dicussed for instance in [7] p.214. The formula for the L-function hereafter follows by (2.1) and (2.2).
\bigskip

{\bf Lemma 2.2:} {\it For $n\equiv m+1\equiv 0 \mod 2$ we get
$$\tau = \bigoplus _{i=1}^{n/2}\bigoplus _{j=1}^{(m-1)/2} (l_i+l'_j,-\kappa ')\oplus (|l_i-l'_j|,-\kappa ')
\oplus \bigoplus _{i=1}^{n/2} (l_i,-\kappa ')$$
hence
$$L(s,\tau ) = \prod _{i,j}\Gamma_{\Bbb C}(s-\kappa '+\frac{l_i+l'_j}{2})\cdot \prod _{i,j;l_i\neq l'_j}\Gamma_{\Bbb C}(s-\kappa '+\frac{|l_i-l'_j|}{2})$$
$$\cdot \prod _{i,j;l_i=l'_j}\Gamma_{\Bbb R}(s-\kappa ') \Gamma_{\Bbb R}(s-\kappa '+1)\cdot \prod _{i=1}^{n/2} \Gamma _{\Bbb C}(s-\kappa '+\frac{l_i}{2}).$$}
\bigskip 

The proof of Lemma 2.2 and  the proof of the following Lemma 2.3 are very much similar to the proof of Lemma 2.1 so we omit it. 

The remaining so-called exceptional case $m\equiv n\equiv 1 \mod 2$ will play a particular role later-on. We put $\epsilon = 0$ or $1$ such that $\epsilon \equiv \delta + \delta ' \mod 2$.
\bigskip

{\bf Lemma 2.3:} {\it In the exceptional case we find
$$\tau = \bigoplus _{i=1}^{(n-1)/2}\bigoplus _{j=1}^{(m-1)/2} (l_i+l'_j,-\kappa ')\oplus (|l_i-l'_j|,-\kappa ')$$
$$\oplus \bigoplus _{i=1}^{(n-1)/2}(l_i,-\kappa ')\oplus \bigoplus _{j=1}^{(m-1)/2}(l'_j,-\kappa ') \oplus (sgn^\epsilon,-\kappa '),$$
hence
$$L(s,\tau )= \prod _{i,j} \Gamma_{\Bbb C}(s-\kappa '+\frac{l_i+l'_j}{2})\cdot \prod _{i,j;l_i\neq l'_j}\Gamma _{\Bbb C}(s-\kappa '+\frac{|l_i-l'_j|}{2})$$
$$\cdot \prod _{i,j;l_i=l'_j}\Gamma _{\Bbb R}(s-\kappa ')\Gamma _{\Bbb R}(s-\kappa '+1)$$
$$\cdot \prod _{i=1}^{(n-1)/2}\Gamma _{\Bbb C}(s-\kappa '+\frac{l_i}{2})\cdot \prod _{j=1}^{(m-1)/2}\Gamma_{\Bbb C}(s-\kappa '+\frac{l'_j}{2}) \cdot \Gamma _{\Bbb R}(s-\kappa '+\epsilon).$$}
\bigskip

{\bf Remark 2.1:} {\it In all cases the L-function $L(s,\check{\tau})$ of the contragredient representation is given by the same formulas as for $L(s,\tau )$ but with $\kappa '$ replaced by $-\kappa '$,
since $\check{\pi}_\infty \cong J(w,l)\otimes sgn^\delta$ and $\check{\sigma}_\infty \cong J(w',l')\otimes sgn^{\delta '}$.}
\bigskip

Next we want to get rid of the unpleasant $\Gamma$-factors attached to the pairs $(i,j)$ where $l_i=l_j'$.
\bigskip

{\bf Lemma 2.4:} {\it If there is a pair $(i,j)$ such that $l_i=l_j'\neq 0$, then for all $t\in \frac{n+m}{2}+\Bbb Z$ 

either $t$ is a pole of $ \Gamma_{\Bbb R}(s-\kappa ') \Gamma _{\Bbb R}(s-\kappa '+1)$

or $t$  is a pole of $\Gamma _{\Bbb R}(1-s+\kappa ') \Gamma _{\Bbb R}(1-s+\kappa '+1)$,

\noindent hence $t$ is a pole of $L(s,\tau )$ or of $L(1-s, \check{\tau})$. In particular $\mathrm{Crit}$ is empty in this case.}
\bigskip

{\bf Proof of Lemma 2.4.} The set of poles of the first product is $\kappa '+\Bbb Z_{\leq 0}$ and that of the second product is $\kappa '+\Bbb N$,
hence each number of the form $\kappa '+z$ with an integer $z$ occurs as a pole. Now the assumption $l_i=l'_j$ implies by means of the congruence 
properties $w+l_i\equiv n+1 \mod 2$ and $w'+l_j'\equiv m+1 \mod 2$ that $w+w'\equiv n+m \mod 2$ hence $\kappa '+\Bbb Z = \frac{n+m}{2} + \Bbb Z$ and therefore any $t$ in this set is a pole
which completes the proof.

Now once we assume the existence of a critical number $t\in \frac{n+m}{2} + \Bbb Z$ the unpleasent $\Gamma$-factors of $L(s,\tau )$ attached to pairs $(i,j)$ with $l_i=l_j'(\neq 0)$ in
the Lemmas 2.1,2,3 cannot turn up by Lemma 2.4. The regularity condition for the remaining explicit $\Gamma$-factors of $L(s,\tau )$ and $L(1-s,\check{\tau})$ says
$|t-\kappa |< L_1+1/2 $ for
$$L_1:= \frac{1}{2} \min _{I=1}^{n/2} \min _{j=1}^{m/2} \{l_i+l_j',|l_i-l_j'|\}$$
in the case $n\equiv m\equiv 0 \mod 2$,
$|t-\kappa |<L_2 +1/2$ for
$$L_2 := \frac{1}{2} \min _{i=1}^{n/2} \min _{j=1}^{(m-1)/2} \{l_i+l_j', |l_i-l_j'|, l_i, l_j'\}$$
in the case $n\equiv m+1\equiv 0 \mod 2$, and eventually $|t-\kappa |<L_3 +1/2$ for
$$L_3 := \frac{1}{2} \min _{i=1}^{(n-1)/2} \min _{j=1}^{(m-1)/2} \{l_i+l_j', |l_i-l_j'|, l_i, l_j' \}$$
in the exceptional case $n\equiv m\equiv 1 \mod 2$, where in addition taking in account the last factor $\Gamma _{\Bbb R}(s-\kappa '+\epsilon)$ of $L(s,\tau )$ we encounter the 
extra condition that $t$ be not a pole of $\Gamma _{\Bbb R}(s-\kappa '+\epsilon) \Gamma _{\Bbb R}(1-s+\kappa '+\epsilon)$. Since we know the poles of $\Gamma (s)$ this is equivalent
to the parity condition 
$$t-\kappa ' \in (2\Bbb N -\epsilon)\cup -(2\Bbb N -1-\epsilon) \leqno(PC_\epsilon)$$
in Proposition 2.1.

Eventually we simplify and unify the bounds $L_i$ (i=1,2,3) by showing that they are in fact all equal to the bound $L$ in the proposition. This is obvious for $L_1$ since
$$|\pm l_i\pm l_j'| \leq l_i+l_j'$$
for all positive $l_i,l_j'$ and hence also for 
$$L_2 = \frac{1}{2} \min _{i=1}^n \min _{j=1}^m \{|l_i-l_j'|\} = L,$$
where for $j=(m+1)/2$ we have $|l_i-l_j'| = |l_i|$. The same observation applies to $L_3 = L$.

In the opposite direction of the proposition for a number $t \in \frac{n+m}{2} +\Bbb Z$ we have $|t-\kappa |\geq 1/2$ or $t=\kappa$, hence $(I)$ implies $L>0$ in the first case. 
In the second case $\frac{n+m}{2} +\kappa$ is integral, i.e. $n+m\equiv w+w'+1 \mod 2$, hence by the congruences $w+l_i\equiv n+1 \mod 2$ and $w'+l_j'\equiv m+1 \mod 2$ we get $l_i+l_j'\equiv n+m+w+w'\equiv 1 \mod 2$ 
and in particular $l_i\neq l_j'$ for all $i,j$, so $L>0$ in this case too. Therefore the L-function has no unpleasant $\Gamma$-factors like those in Lemma 2.4 and $t$ satisfying $(I)$ and $(PC)_\epsilon$ in addition
in the exceptional case is a critical number. This completes the proof of the proposition.
\bigskip 

{\bf Remark 2.2:} {\it There is an obvious reflection map 
$$\theta : \mathrm{Crit} \longrightarrow \mathrm{Crit}, t \mapsto w+w' + 1 -t.$$}
\bigskip

{\bf Proposition 2.2:} {\it If the set $\mathrm{Crit}$ is non empty then for $i=1,...,n $ and $j=1,...,m$ 
we have for each pair $(i,j)$ the alternative 
$$l_i \neq l'_j \mbox{ or } l_i = l'_j = 0.$$
The second case only occurs for $n\equiv m \equiv 1 \mod 2$  with $i = \frac{n+1}{2}, j = \frac{m+1}{2}$.}

\bigskip 

{\bf Proof.} The existence of a critical $t$ implies $L \neq 0$ hence in particular $L_0 \neq 0$, so 
$l_i \neq l'_j$ for $i \leq \frac{n}{2}$ and $j \leq \frac{m}{2}$ which are the positive components of
$l$ and $l'$. Therefore certainly we have $l_i \neq -l'_j = l'_{m+1-j}$ and $l'_j \neq -l_i = l_{n+1-i}$. So 
the only remaining case is $l_i = l'_j = 0$ as described.

\bigskip

In the non-exceptional case we can reverse the conclusion in proposition 2.2.

\bigskip

{\bf Proposition 2.3:} {\it If we are not in the exceptional case $n\equiv m\equiv 1\mod 2$ the set of critical numbers is
non empty if and only if for each pair $(i,j)$ we have $l_i \neq l'_j$, i.e. $L_0 \neq 0$. In this situation the two numbers
$t_0 := L -1 + \frac{w+w'+1}{2}$ and $\theta (t_0)$ are critical.}

\bigskip

{\bf Proof.} We are left to show the existence of a critical number if $L_0$ does not vanish. Now for $L_0\neq 0$ obviously
$t_0$ is critical.

From now on we will assume the existence of a critical number. Since we excluded the case $n = m = 1$ in this paper we may and will assume that the Langlands parameters satisfy the
\bigskip

{\bf Hypothesis:} $l_1 > l'_1$.
\bigskip 

The correlation of the Langlands parameters $l$ and $l'$  is recorded by the {\it position tupel} of integers
$a = (a_1,...,a_m)$ uniquely determined by the property
$$l_{a_j} > l'_j \geq l_{1+a_j} \mbox{ for } j=1,...,m.$$
Obviously we have
$$1 \leq a_1 \leq a_2 \leq ... \leq a_m \leq n-1.$$
In general the sequence $a_j$ is not constant.
\bigskip

{\bf Lemma 2.5:} {\it The position tupel $a$  enjoys the symmetry
$$ a_j + a_{m+1-j} = n \mbox{ for } j=1,...,m \leqno(2.3)$$
except for $j = \frac{m+1}{2}$ in the exceptional case $m \equiv n \equiv 1 \mod 2$ where $a_j = \frac{n-1}{2}$. }
\bigskip  

This follows immediately from the defining inequalities via the symmetry of the Langlands parameters 
$$l_i + l_{n+1-i} = 0 \mbox{ and } l'_j + l'_{m+1-j} = 0$$
\bigskip

{\bf Remark 2.3:} {\it The constant case: If $m>1$ and $a_1 = a_m$ necessarily $n$ must be even, $a_i = n/2$ and
$$l_{\frac{n}{2}} > l'_1 > l'_2 > ... > l'_m > -l_{\frac{n}{2}}.$$
Assuming that $n$ is even we arrive at the same conclusions for $m=1$. For odd $n$ and $m=1$ we get $a_1 = \frac{n-1}{2}$.} 
\bigskip

If $a$ is not constant there are {\it jump indices}. We call $j$ a {\it jump index} if we have  $a_j < a_{j+1}$. Let $\{j_1,...,j_k\}$ denote the set of {\it jump indices} and
let $k:=0$ for constant $a$. Note that for $m=1$ certainly $a$ is constant hence $k=0$. Suppose we have $m>1$. Then
by the symmetry property (2.3) we get 
\bigskip

{\bf Remark 2.4:} {\it In the non-exceptional case $j$ is a jump index if and only if $m-j$ is a jump index. In the exceptional case this holds for $j\neq \frac{m\pm 1}{2}$ and moreover 
$j = \frac{m+1}{2}$ is always a jump index whereas $\frac{m-1}{2}$ is a jump index if and only if we have $l_{\frac{n-1}{2}} < l'_{\frac{m-1}{2}}$, which therefore is equivalent to $k$ being even. Furthermore 
except for the exceptional case with odd $k$ the jump indices 
enjoy the symmetry property
$$j_\kappa + j_{k+1-\kappa} = m \mbox{ for } \kappa  = 1,...,k. \leqno(2.4)$$
In the excluded case this symmetry still holds for $\kappa \neq \frac{k+1}{2}$ whereas we have $j_\frac{k+1}{2} = \frac{m+1}{2}$.}

\bigskip

In terms of the position tupel $a$ and  we can  express in the non-exceptional case the bound $L$ in Proposition 2.1 as
$$L = \frac{1}{2} \min \{ l_{a_j} - l'_j, l'_j - l_{1+a_j}; j=1,...,m\}$$
hence the inequality $(I)$ in Proposition 2.1  for (and so the characterization of) critical numbers transforms into the simultaneous inequalities
$$\pm (t - \frac{w+w'-1}{2}) < \frac{1}{2}(l_{a_j} -l'_j + 1), \frac{1}{2}(l'_j - l_{1+a_j} + 1)$$
for $j=1,...,m$. In the exceptional case we get
$$L = \frac{1}{2} \min \{l'_\frac{m-1}{2}, l_\frac{n-1}{2}, l_{a_j} - l'_j, l'_j - l_{1+a_j}; j\neq \frac{m+1}{2} \}$$
and the respective modification of the inequalities for critical numbers.

\section{Heighest weights and critical numbers}

For arbitrary $n$ we consider irreducible algebraic representations $(\rho _\mu ,M_\mu)$ of $GL_n$ over $\Bbb Q$ with dominant highest weight $\mu$, i.e. $\mu = (\mu _1,...,\mu _n)$ in
$$X^+(n) := \{\mu \in \Bbb Z ^n; \mu _1\geq \mu _2\geq ...\geq \mu _n \}.$$
We are in particular interested in the essentially selfdual representations. These are precisely those representations which satisfy the {\it purity condition} saying that for a suitable integer 
$w = wt(\mu )$ the so called {\it weight} of $\mu$ we have
$$\mu _i + \mu_{n+1-i} = w \mbox{ for all } i.$$
Recall that by (3.5) and (3.6) in [8] the set $X^+_0(n)$ of these $\mu$ corresponds to the set of Langlands parameters $L^+_0(n)$ via the bijection 
$$L^+_0(n) \longrightarrow X^+_0(n)$$
$$(w,l) \mapsto \mu = (\frac{w+l_1+1-n}{2},\frac{w+l_2+3-n}{2},...,\frac{w+l_n+n-1}{2}),$$
where the inverse mapping is given by $\mu \mapsto (w,l)$ with $w = \mu _1 + \mu _n$ and $l_i = 2 \mu _i + n + 1 - w -2i$, and where any $\pi \in Coh(\mu )$ has infinity component
$$\pi _\infty \cong J(-w,l)\otimes sgn^\delta$$ for suitable $\delta$. Via this bijection we want to express
the previous simultaneous inequalities for critical numbers in terms of the associated highest weights. So let $\mu \in X_0^+(n)$ and $\nu \in X_0^+(m)$ correspond to the Langlands parameters 
$(w,l)$ and $(w',l')$ as considered in the previous section for $(\pi ,\sigma )\in \Coh(\mu )\times \Coh(\nu )$ with a critical pair $\pi _\infty ,\sigma _\infty )$. Using the position tupel $a$ we introduce the weight vector $\lambda = (\lambda _1,...,\lambda _m)$ where we put
$$\lambda _j := \nu _j + a_j - j \mbox{ for } j=1,...,m.$$
\bigskip

{\bf Proposition 3.1:} {\it The weight $\lambda$ is dominant and satisfies the purity condition 
$$\lambda _j + \lambda _{m+1-j} = w' + n - m -1 \leqno(P_j)$$
for $j \neq \frac{m+1}{2}$. In the exceptional case $n\equiv m\equiv 1\mod 2$ we have 
$$2 \lambda _\frac{m+1}{2} = w' + n - m - 2,$$
whereas in the non-exceptional case $(P_j)$ holds for all $j$.}

\bigskip

An obvious consequence is

\bigskip

{\bf Corollary 3.1:} {\it In the non-exceptional case we have $\lambda \in X_0^+(m)$. In the exceptional case for $m>1$ we may modify $\lambda$ to achieve the purity property for
$$\lambda _{mod} := (\lambda _1,...,\lambda _\frac{m+1}{2},\lambda _\frac{m+3}{2} - 1,...,\lambda _m - 1) \in X_0^+(m)$$
and the truncated
$$\lambda _{tr} := (\lambda _1,...,\lambda _\frac{m-1}{2},\lambda _\frac{m+3}{2},...,\lambda _m) \in X_0^+(m-1).$$
For $m=1$ always $\lambda \in X_0^+(1)$.}

\bigskip

{\bf Proof.} If $j$ is not a jump index we have $a_j=a_{j+1}$ hence
$$\lambda _j = \nu _j + a_j - j \geq \nu _{j+1} + a_{j+1} - j > \lambda _{j+1}$$
since $\nu$ is dominant. Otherwise we have $a_j + 1 \leq a_{j+1}$ so the inequalities defining the position tupel $a$ in particular imply that we have
$l'_j \geq l_{1+a_j}$ and $l_{a_{j+1}} > l'_{j+1}$. In terms of the associated highest weights this says
$$2\nu _j + m + 1 - w' -2j \geq 2\mu _{1+a_j} + n + 1 - w - 2(1+a_j)$$
and  
$$2\mu _{a_{j+1}} + n + 1 - w - 2a_{j+1} > 2\nu _{j+1} + m + 1 - w' - 2(j+1).$$
So with $R:= \frac{1}{2}(m-n+w-w')$ we get 
$$\lambda _j + R \geq \mu _{1+a_j} - 1  \mbox{ and } \mu _{a_{j+1}} > \lambda _{j+1} + R.$$ 
Since 
$\mu _{1+a_j} \geq \mu _{a_{j+1}}$ this implies $\lambda _j \geq \lambda _{j+1}$ which proves that $\lambda$ is dominant.

In the non-exceptional case the purity property $(P_j)$ follows directly from the purity of $\nu$ and the symmetry property $(AS)$ of the position tupel $a$.
In the exceptional case the same argument works for $j\neq \frac{m+1}{2}$ and for $j=\frac{m+1}{2}$ the purity of $\nu$  and the identity $a_\frac{m+1}{2} = \frac{n-1}{2}$ 
supply the claimed extra formula for $\lambda _\frac{m+1}{2}$.

\bigskip

We now turn to the {\it dual} highest weight $\check{\mu}$ with the components $\check{\mu } _j = -\mu _{n+1-j}$ and corresponding  Langlands parameter $(-w,l)$.

\bigskip

{\bf Proposition 3.2:} {\it In the non-exceptional case a number $t \in \frac{m+n}{2} + \Bbb Z$ satisfies the inequality $(I)$ hence is critical if and only if the integer
$s := t + \frac{n-m}{2} - 1 $ satifies the simultaneous inequalities 
$$\check{\mu }_{a_j} \geq \lambda _j - s \geq \check{\mu }_{1+a_j} \mbox{ for } j=1,...,m. \leqno(I_j)$$
In the exceptional case the condition is that $s$ fulfills $(I_j)$ for all $j \neq \frac{m+1}{2}$ and
the ``middle condition'' 
$$\check{\mu}_\frac{n-1}{2} \geq \lambda _\frac{m+1}{2} - s \geq \check{\mu}_\frac{n+3}{2} -1 \mbox{ for even }k,$$
respectively
$$\lambda _\frac{m-1}{2} \geq \check{\mu}_\frac{n+1}{2} + s \geq \lambda _\frac{m+3}{2} \mbox{ for odd }k.$$}
\bigskip

{\bf Proof.} As we observed at the end of the previous section a number $t$ satisfying $(I)$ is characterized by a system of inequalities
envolving certain Langlands parameters  $l_i$ and $l'_j$ which we now express in terms of the associated weights
$$l_i = 2\check{\mu}_i + n + 1 + w - 2i$$
respectively
$$l'_j = 2\nu _j + m + 1 - w' -2j.$$
In the non-exceptional case this system gets the form
$$\pm (t - \frac{w+w'+1}{2}) < \check{\mu}_{a_j} - \nu _j - a_j + j + \frac{n-m}{2} + \frac{w+w'+1}{2} $$
and
$$\pm (t - \frac{w+w'+1}{2}) < \nu _j + a_j - j - \check{\mu}_{1+a_j} - \frac{n-m}{2} - \frac{w+w'+1}{2} + 2$$
for $j=1,...,m$. In the exceptional case the system consists of these inequalities for $j\neq \frac{m+1}{2}$ and the inequality
$$\pm (t - \frac{w+w'+1}{2}) < \check{\mu }_\frac{n-1}{2} + \frac{3+w}{2} \mbox{ for even }k$$
respectively
$$\pm (t - \frac{w+w'+1}{2}) < \nu _\frac{m-1}{2} + \frac{3-w'}{2} \mbox{ for odd }k.$$
By the purity property of $\nu$ and $\check{\mu}$ our system of inequalities reduces to the simultaneous inequalities $(I_j)$ and the last 
inequality transforms into the ``middle condition'' which completes the proof.

\bigskip

{\bf Remark 3.1:} {\it In the exceptional case for odd $k$ the inequalities $(I_\frac{m-1}{2})$ and $(I_\frac{m+3}{2})$ imply the ``middle condition'' 
hence the latter can be deleted in the criterion of the proposition.}

\bigskip

Taking into account the jump indices $j_\kappa$ of the position tupel $a$ we can easily reformulate Proposition 3.2. Set $j_0 := 0$ and $j_{k+1} := m$.

\bigskip

{\bf Proposition 3.3:} {\it  In the non-exceptional case a number $t \in \frac{m+n}{2} + \Bbb Z$ satisfies the inequality $(I)$ (hence is critical) if and only if the integer
$s := t + \frac{n-m}{2} - 1 $ satifies the simultaneous inequalities 
$$\check{\mu}_{a_{j_\kappa}} \geq \lambda _{1+ j_{\kappa -1}} -s \geq ...\geq \lambda _{j_\kappa} -s \geq \check{\mu}_{1+{a_{j_\kappa}}} \leqno(S_\kappa)$$
for $\kappa = 1,...,k+1$. Furthermore the inequalities of the listed $\lambda$-components are strict, i.e.
$$\lambda _{1+j_{\kappa -1}} >  \lambda _{2+j_{\kappa -1}} > ... > \lambda _{j_\kappa} .$$
In the exceptional case $(I)$ holds for odd $k$ if and only if the integer $s$  satisfies $(S_\kappa)$ for all $\kappa \neq \frac{k+1}{2}$ and for $\kappa = \frac{k+1}{2}$ the truncated series of inequalities
$$\check{\mu}_\frac{n-1}{2} \geq \lambda _{1+j_\frac{k-1}{2}} -s \geq ... \geq \lambda _{j_\frac{k+1}{2}-1} -s \geq \check{\mu}_\frac{n+1}{2} \leqno(S'_\frac{k+1}{2})$$
where $j_\frac{k+1}{2} = \frac{m+1}{2}$ and $a_\frac{m+1}{2} = a_\frac{m-1}{2} = \frac{n-1}{2}$. For even $k$ the integer $s$ has to satisfy $(S_\kappa)$ for all $\kappa \neq 1+\frac{k}{2}$ and the 
``middle condition''
$$\check{\mu}_{\frac{n-1}{2}} \geq \lambda _{\frac{m+1}{2}} - s \geq \check{\mu}_{\frac{n+3}{2}} -1.$$ }

\bigskip

\section{Inequality systems and ramification}

Our goal in this section is to interprete via the ramification law of restricted irreducible representations the simultaneous inequalities in Proposition 3.3 as embeddings of representations attached to 
suitable choices of highest weights. To begin with we observe that the essential part of the inequality sequence  $(S_\kappa)$ is the equivalent four term inequality
$$\check{\mu}_{a_{j_\kappa}} \geq \lambda _{1+ j_{\kappa -1}} -s \geq \lambda _{j_\kappa} -s \geq \check{\mu}_{1+{a_{j_\kappa}}}.\leqno(4.1)$$
With that in mind we are led to study arbitrary systems of inequalities  of the form
$$u_{2\rho -1} \geq v_{2\rho -1} \geq v_{2\rho} \geq u_{2\rho} \mbox{ for } \rho =1,...,r \leqno(4.2)$$
where $u=(u_1,...,u_{2r})$ and $v=(v_1,...,v_{2r})$ are dominant weights for $GL_{2r}$ with purity property, i.e. $u,v \in X_0^+(2r)$. In this situation we say that the pair $(u,v)$ of heighest weights 
satisfies $(4.2)$.
In the total $4r$-term inequality sequence $(4.2)$ we must take care of the {\it middle group} 
$$v_r \geq u_r \geq u_{r+1} \geq v_{r+1} \mbox{ \ \ for even \ \ } r,$$
respectively
$$u_r \geq v_r \geq v_{r+1} \geq u_{r+1}  \mbox{ \ \ for odd \ \ } r.$$
For a highest weight $y \in X_0^+(2r)$ we call the difference of the middle components $d(y) := y_r - y_{r+1}$ the {\it defect} of $y$.
\bigskip  

{\bf Remark 4.1:} {\it For a pair $(u,v)$ of $u,v \in X_0^+(2r)$ satisfying $(4.2)$ we necessarily have $d(u) \leq d(v)$ for even $r$ and $d(v) \leq d(u)$ for odd $r$.   }
\bigskip

For technical reasons we suppose for the time being that the middle group is {\it special} in the sense that $d(u)=0$ or $d(v)=0$. For the set of such special heighest weights we 
introduce the notation
$$X_{sp}(2r) := \{y \in X_0^+(2r); d(y)=0\}.$$
We will consider ``splitting maps'' $\theta =(\theta _1,\theta _2), \theta '=(\theta _1',\theta _2')$ with 
$$\theta : X_{sp}(2r) \longrightarrow X_0^+(r+1) \times X_0^+(r+1),$$
$$\theta ': X_0^+(2r) \longrightarrow X_0^+(r) \times X_0^+(r)$$
for even $r$ and
$$\theta : X_0^+(2r) \longrightarrow X_0^+(r+1) \times X_0^+(r+1),$$
$$\theta ':= X_{sp}(2r) \longrightarrow X_0^+(r) \times X_0^+(r)$$
for odd $r$ defined as follows. For $y \in X_{sp}(2r)$ and $z \in X_0^+(2r)$ we put for even $r$
$$\theta _1(y) := (y_1,y_2,y_4,...,y_r=y_{r+1},y_{r+3},...,y_{2r-1},y_{2r}),$$
$$\theta _2(y) := (y_1,y_3,...,y_{r-1},y_{r+1}=y_{r},y_{r+2},y_{r+4},...,y_{2r})$$
and 
$$\theta _1'(z) := (z_2,z_4,...,z_r,z_{r+1},z_{r+3},...,z_{2r-1}),$$
$$\theta _2'(z) := (z_1,z_3,...,z_{r-1},z_{r+2},z_{r+4},...,z_{2r}).$$
For odd $r$ we similarly define 
$$\theta _1(z) := (z_1,z_2,z_4,...,z_{r-1},z_{r+2},z_{r+4},...,z_{2r-1},z_{2r}),$$
$$\theta _2(z) := (z_1,z_3,...,z_r,z_{r+1},z_{r+3},...,z_{2r})$$
and
$$\theta _1'(y) := (y_2,y_4,...,y_{r+1}=y_{r},y_{r+2},y_{r+4},...,y_{2r-1}),$$
$$\theta _2'(y) := (y_1,y_3,...,y_r=y_{r+1},y_{r+3},...,y_{2r}).$$
In particular for $r=1$ where $y_1=y_2$ we take $\theta _i(z)=z$ and $\theta _i'(y)=y_1$ for $i=1,2$.
\bigskip

{\bf Remark 4.2:} {\it Obviously the images of $\theta$ and $\theta'$ consist of dominant weights and they inherit the purity property of $y$ and $z$.  }
\bigskip

Our main technical tool later-on will be
\bigskip

{\bf Proposition 4.1:} {\it For a given pair $(y,z)$ of highest weights for $GL_{2r}$ with $y\in X_{sp}(2r)$ and $z\in X_0^+(2r)$ the following statements are equivalent:

a) The pair $(y,z)$ satisfies $(4.2)$ for even $r$ respectively the pair $(z,y)$ satisfies $(4.2)$ for odd $r$.

b) The images under the splitting maps $\theta$ and $\theta '$ have the embedding property
$$\theta _i'(z) \prec \theta _i(y) \ \ (i=1,2) \mbox{ \ \ for even \ \ } r$$
respectively
$$\theta _i'(y) \prec \theta _i(z) \ \  (i=1,2) \mbox{ \ \ for odd \ \ } r.$$}
\bigskip

{\bf Proof.} We first work out the proof for even $r$. The inequality system (4.2) for the pair $(y,z)$ may be separated into two groups of
inequalities. The first group is 
$$(y_1 \geq ) z_2 \geq y_2,$$
$$(y_i \geq ) z_{i+2} \geq y_{i+2} \mbox{ \ \ for even \  } i, \ 2 \leq i \leq r-2,$$
$$y_i \geq z_i (\geq y_{i+2} )  \mbox{ \ \ for odd \  } i, \ r+1\leq i \leq 2r-3,$$
$$y_{2r-1} \geq z_{2r-1} (\geq y_{2r}).$$
The second group is 
$$y_i \geq z_i (\geq y_{i+2})  \mbox{ \ \ for odd \  } i, \ 1\leq i \leq r-1,$$
$$(y_{i-2} \geq ) z_i \geq y_i  \mbox{ \ \ for even \  } i, \ r+2 \leq i \leq 2r.$$
Now the first group precisely describes the embedding statement $\theta _1'(z) \prec \theta _1(y)$ and the second group expresses the statement $\theta _2'(z) \prec \theta _2(y)$,
hence the statements a) and b) are equivalent.

For odd $r$ the inequality system $(4.2)$ for the pair $(z,y)$ similarly seperates into two groups. Now the first group is
$$(z_1\geq ) y_2 \geq z_2,$$
$$(z_{i-2}\geq ) y_i \geq z_i \mbox{ \ \ for even \  }i,\ 4\leq i\leq r-1,$$
$$z_i \geq y_i (\geq z_{i+2}) \mbox{ \ \ for odd \ } i, \ r+2\leq i\leq 2r-3,$$
$$z_{2r-1} \geq y_{2r-1} (\geq z_{2r}).$$
The second group is
$$z_i \geq y_i (\geq z_{i+2}) \mbox{ \ \ for odd \  } i, \ 1\leq i \leq r,$$
$$(z_{i-2} \geq ) y_i \geq z_i \mbox{ \ \ for even \  } i, \ r+1\leq i \leq 2r.$$
Here the first group describes the embedding statement $\theta _1'(y) \prec \theta _1(z)$ whereas the second group says $\theta _2'(y) \prec \theta _2(z)$, hence like in the even case
the statements a) and b) are equivalent and the proof is complete.

We can easily extend the criterion of Proposition 4.1 to an arbitrary pair $(u,v)$ of heighest weights $u,v\in X_0^+(2r)$. Let $l:=(0,...,0,1,...,1)\in \Bbb Z ^{2r}$ have a zero in each of the first $r$
components and $1$ in each of the remaining $r$ components. Further let $d$ denote the minimum of the defects $d(u),d(v)$, i.e.
$$d:=min\{ u_r-u_{r+1}, v_r - v_{r+1} \},$$
and define the modified weights $\hat{u} := u + d\cdot l$ and $\hat{v} := v + d \cdot l$ such that in particular we have $d(\hat{u})=0$ or $d(\hat{v})=0$.
\bigskip

{\bf Corollary 4.2:} {\it The pair $(u,v)$ satisfies $(4.2)$ if and only if $d(\hat{u})=0$ for even $r$ respectively $d(\hat{v})=0$ for odd $r$ and if we have the two embeddings 
$$\theta _i'(\hat{v}) \prec \theta _i(\hat{u})  \mbox{ \ \ for \ \ } i=1,2.$$} 
\bigskip

{\bf Proof.} Certainly $\hat{u}$ and $\hat{v}$ still belong to $X_0^+(2r)$ and $(u,v)$ satisfies $(4.2)$ if and only if $(\hat{u},\hat{v})$ does. Now we can apply  Proposition 4.1 to $(\hat{u},\hat{v})$.
If $(4.2)$ holds for $(\hat{u},\hat{v})$ then necessarily $d(\hat{u})=0$ for even $r$ and $d(\hat{v})=0$ for odd $r$ by Remark 4.1. In both cases we get $\theta _i'(\hat{v}) \prec \theta _i(\hat{u})$ for 
$i=1,2$ by Proposition 4.1. Conversely if we have $d(\hat{u})=0$ for even $r$ respectively $d(\hat{v})=0$ for odd $r$ and $\theta _i'(\hat{v}) \prec \theta _i(\hat{u})$ for $i=1,2$ then the proposition yields
$(4.2)$ for $(\hat{u},\hat{v})$ in the even case and in the odd case as well.

\section{Critical numbers and embeddings}

We now want to apply the technique from the previous section in order to parametrize the set $\mathrm{Crit}$ of critical numbers of a critical pair $(\pi _\infty ,\sigma _\infty) $ attached to
$(\pi ,\sigma ) \in \Coh(\mu )\times \Coh(\nu )$. 
For this purpose we will use the system of simultaneous inequalities $(S_\kappa)$ in Proposition 3.3 which we must reformulate slightly modified in a suitable form adjusted to the
requirements of Corollary 4.2..

The non-exceptional case being less involved we first deal with this case. In view of Proposition 3.3 and $(4.1)$ we  put 
$$u_{2\rho -1} := \check{\mu}_{a_{j_\rho}}, \ \ u_{2\rho} := \check{\mu}_{1+a_{j_\rho}}, $$
and
$$v_{2\rho -1} := \lambda _{1+j_{\rho -1}} - s, \ \ v_{2\rho} := \lambda _{j_\rho} - s $$
for $\rho =1,...,r$. Since $v$ depends on $s$ we also write $v(s)=v$.
\bigskip

{\bf Lemma 5.1:} {\it The $r$-tupel $u$ and $v$ are highest weights with purity property , i.e. $u,v\in X_0^+(2r)$ of respective weight $wt(u)=-w$ and $wt(v)=w'+n-m-1-2s$.   }
\bigskip

{\bf Proof.} Obviously $u$ and $v$ are dominant weights in $X^+(2r)$ since $\check{\mu}$ like $\mu$ is dominant, $\lambda$ is dominant by Proposition 3.1 and we have
$$a_{j_1} < a_{j_2} < ... <a_{j_r}$$
by the definition of jump indices. For the purity of $u$ we have to show
$$u_{2\rho -1} + u_{2(r+1-\rho )} = -w   \mbox{ \ \ for \ \ } \rho =1,...,r.$$
This is an easy consequence of the purity of $\check{\mu}$, since combining the symmetry $(2.4)$ of jump indices with the symmetry $(2.3)$ of the position tupel $a$ 
we get (remember $j_0=0$ and $j_r=m$)
$$a_{j_\rho} + a_{j_{1+r-\rho}} = n \mbox{\ \ for \ \ }\rho =1,...,r, \leqno(5.1)$$
where we use $a_{1+j_\rho} = a_{j_{\rho + 1}}$. For the purity of $v$ we exploit the purity of $\lambda \in X_0^+(m)$ to show
$$v_{2\rho -1} + v_{2(r+1-\rho )} = w' + n - m - 1 -2s \mbox{ \ \ for \ \ }\rho =1,...,r.\leqno(5.2)$$
Again the symmetry $(2.4)$ of jump indices and Proposition 3.1 imply
$$v_{2\rho -1} + v_{2(r+1-\rho )} = \lambda _{1+j_{\rho -1}} - s + \lambda _{j_{r+1-\rho}} - s = w'+n-m-1-2s,$$
as required.
\bigskip

{\bf Remark 5.1:} {\it The respective defects of $u$ and $v$ are
$$d(u) = \check{\mu}_{1+a_{j_\frac{r}{2}}} - \check{\mu}_{a_{j_{1+\frac{r}{2}}}}  \mbox{ \ \ for even\ \ }r$$
and
$$d(v)=\lambda _{1+j_\frac{r-1}{2}} - \lambda _{j_\frac{r+1}{2}}  \mbox{ \ \ for odd \ \ } r.$$
In particular $d(v)$ is independent of $s$ and the existence of a critical number implies $d(u)\leq d(v)$ for even $r$ and $d(v)\leq d(u)$ for odd $r$.}
\bigskip

{\bf Proof.} The explicit form of the defect can directly be read off from the definition and the inequality follows by $(4.2)$ for the middle group.

Now with $d:=\min\{d(u), d(v)\}$ we define the modified $\hat{u} := u +  d\cdot l$ and $\hat{v} := v + d\cdot l$ like in the previous section. Passing to the dual of $\hat{u}$ we find 
$\check{\hat{u}} = \hat{u} + w - d$.
\bigskip

{\bf Remark 5.2:} {\it For the dual $\tilde{\mu} := \check{\hat{u}}$ of $\hat{u}$ we explicitely have
$$\tilde{\mu}_{2\rho -1} = \left\{\begin{array}{ll}
                                   \mu _{a_{j_\rho}} - d  & \mbox{for \ $2\rho -1 \leq r,$}\\
                                   \mu _{a_{j_\rho}}      & \mbox{for \ $2\rho -1 > r,$}\\
                                  \end{array}
                                  \right.$$
$$\tilde{\mu}_{2\rho} = \left\{\begin{array}{ll}
                                \mu _{1+a_{j_\rho}} - d   & \mbox{for \ $2\rho \leq r,$}\\
                                \mu _{1+a_{j_\rho}}       & \mbox{for \ $2\rho > r.$}\\
                               \end{array}
                                  \right.$$}
\bigskip

{\bf Proof.} By Lemma 5.1 we know $wt(u)=-w$ hence $\hat{u}$ has weight $wt(\hat{u}) = d-w$ and the dual of $\hat{u}$ is  $\tilde{\mu} = \check{\hat{u}} = \hat{u} + w - d$. By definition
of $\hat{u}$ we have
$$\hat{u}_{2\rho -1} = \left\{\begin{array}{ll}
                               \check{\mu}_{a_{j_\rho}}      & \mbox{for \ $2\rho -1 \leq r,$}\\
                               \check{\mu}_{a_{j_\rho}} + d  & \mbox{for \ $2\rho - 1 > r,$}\\
                              \end{array}
                              \right.$$
and 
$$\hat{u}_{2\rho} = \left\{\begin{array}{ll}
                            \check{\mu}_{1+a_{j_\rho}}       & \mbox{for \ $2\rho \leq r,$}\\
                            \check{\mu}_{1+a_{j_\rho}} + d   & \mbox{for \ $2\rho > r.$}\\
                           \end{array}
                           \right.$$
Since $\check{\mu} = \mu - w$ we arrive at the claimed formula for $\tilde{\mu}$.

In a similar way by definition of $\hat{v} = v+d\cdot l$ we have
$$\hat{v}_{2\rho -1} = \left\{\begin{array}{ll}
                               \lambda _{1+j_{\rho -1}} - s      & \mbox{for \ $2\rho -1 \leq r,$}\\
                               \lambda _{1+j_{\rho -1}} + d - s  & \mbox{for \ $2\rho -1 > r,$}\\
                              \end{array}
                              \right.$$
and
$$\hat{v}_{2\rho} = \left\{\begin{array}{ll}
                            \lambda _{j_\rho} -s                 & \mbox{for \ $2\rho \leq r,$}\\
                            \lambda _{j_\rho} + d - s            & \mbox{for \ $2\rho > r.$}\\
                           \end{array}
                           \right.$$
Isolating the influence of $s$ in these formulas we put $\tilde{\lambda} := \hat{v} + s$ and find the obvious explicit description of $\tilde{\lambda}$ from the preceeding lines. 
\bigskip
                            
{\bf Theorem A:} {\it In the non-exceptional case the mapping $t \mapsto t + \frac{n-m}{2} - 1$ sets up a bijection 
$$\mathrm{Crit}(\pi_\infty ,\sigma _\infty) \longrightarrow \bigcap _{i=1,2}\Emb(\theta _i'(\tilde{\lambda}),\theta _i(\check{\tilde{\mu}})).$$}
\bigskip

{\bf Proof.} For a critical pair $(\pi _\infty ,\sigma _\infty )$ by Proposition 3.3 and $(4.1)$ a number $t\in \frac{m+n}{2} + \Bbb Z$ is critical if and only if for $s=t+\frac{n-m}{2}-1$ our concrete
$(u,v)$ in the beginning of this section satisfies $(4.2)$. By Corollary 4.2 this is equivalent to the embedding property 
$$\theta _i'(\hat{v}) \ \  \prec  \ \ \theta _i(\hat{u})  \mbox{ \ \  for \ \ } i=1,2.$$
In terms of the previously defined weights $\tilde{\lambda}$ and $\tilde{\mu}$ this says
$$\theta _i'(\tilde{\lambda}) - s  \  \ \prec \ \ \theta _i(\check{\tilde{\mu}}),\leqno(5.3)$$
i.e.
$$s \in  \bigcap _{i=1,2} \Emb(\theta _i'(\tilde{\lambda}), \theta _i(\check{\tilde{\mu}})),$$
which settles the proof of the theorem.
\bigskip

{\bf Corollary 5.1:} {\it In terms of the Tate modules $T_r(s)=\det ^s$ for $GL_r$ we have
$$\bigcap _{i=1,2} (M_{\theta _i'(\tilde{\lambda})} \otimes M_{\theta _i(\tilde{\mu})})^{SL_r} \cong \bigoplus _{t\in Crit} T_r(t + \frac{n-m}{2} - 1).$$}
\bigskip

{\bf Proof.} This follows by the general principle that for $\alpha \in X_0^+(r)$ and $\beta \in X_0^+(r+1)$ we have $\alpha - s \prec \check{\beta}$ if and only if the Tate module $T_r(s)$ is a
($1$-dimensional) $GL_r$-submodule of $(M_\alpha \otimes M_\beta )^{SL_r}$.

For the remainder of this section we suppose that we are in the exceptional case. In view of Proposition 2.1 the eventual parametrization of the critical numbers will be
provided by a parameter set of embeddable twists similar to the non-exceptional case subject to an additional parity condition that takes care of $(PC)_\epsilon$. Again as 
before we must distinguish the two cases where the number $k$ of jump indices of the position tupel $a$ is even or odd (i.e. where $r=k+1$ is odd or even).

{\it Case $r\equiv 1 \mod 2$.} Recall that by Proposition 3.3 the necessary condition $(I)$ for a critical number $t$ is equivalent to the system of inequalities $(S_\kappa )$
for all $\kappa \neq \frac{r+1}{2}$ plus the middle condition as formulated in Proposition 3.3. As before an equivalent formulation is provided by a system of inequalities of the form $(4.2)$
for highest weights $u,v \in X_0^+(2r)$ now given by
$$u_{2\rho -1} := \check{\mu}_{a_{j_\rho}}, \ \ u_{2\rho } := \check{\mu}_{1+a_{j_\rho}} \mbox{ \ \ for \ \ } \rho =1,...,\frac{r-1}{2},$$
$$u_r =: \check{\mu}_{\frac{n-1}{2}}, \ \ u_{r+1} := \check{\mu}_{\frac{n+3}{2}} -1 ,$$
$$u_{2\rho -1} := \check{\mu}_{a_{j_\rho}} - 1, \ \ u_{2\rho} := \check{\mu}_{1+a_{j_\rho}} - 1 \mbox{ \ \ for \ \ } \rho =\frac{r+3}{2},...,r,$$

$$v_{2\rho -1} := \lambda _{1+j_{\rho -1}} - s, \ \ v_{2\rho} := \lambda _{j_\rho} - s \mbox{ \ \ for \ \ } \rho =1,...,\frac{r-1}{2},$$
$$v_r := \lambda _{\frac{m+1}{2}} - s =: v_{r+1},$$
$$v_{2\rho -1} := \lambda _{1+j_{\rho -1}} -1 -s, \ \ v_{2\rho} := \lambda _{j_\rho} -1 -s \mbox{ \ \ for \ \ } \rho =\frac{r+3}{2},...,r.$$

The same arguments as in the proof of Lemma 5.1 together with Proposition 3.1 show
\bigskip

{\bf Lemma 5.2:} {\it The weights $u$ and $v$ belong to $X_0^+(2r)$ and have weight
$$wt(u) = -w-1, \ \ wt(v) = w'+n-m-2-s.$$}
\bigskip

By definition here $v$ has defect $d(v)=0$, i.e. $v\in X_{sp}(2r)$ and therefore by Proposition 4.1 the pair $(u,v)$ satisfies $(4.2)$ if and only if we have the embeddings
$$\theta _i'(v)  \ \ \prec  \ \  \theta _i(u) \mbox{ \ \ for \ \ } i=1,2.$$
Like in the non-exceptional case using the notation $\tilde{\lambda} := v + s$ (which is independent of $s$) and $\tilde{\mu} :=\check{u}$ this again says
$$\theta _i'(\tilde{\lambda}) - s  \ \ \prec  \ \ \theta _i(\check{\tilde{\mu}}) \mbox{ \ \ for \ \ } i=1,2. \leqno(5.3)$$
\bigskip

{\it Case $r\equiv 0 \mod 2$.} Here condition $(I)$ again by Proposition 3.3 is equivalent to the system of inequalities $(S_\kappa )$ for all $\kappa \neq \frac{r}{2}$
plus the truncated series of inequalities $(S'_\frac{r}{2})$. Thus we are led to define

$$u_{2\rho -1} := \check{\mu}_{a_{j_\rho}}, \ \ u_{2\rho} := \check{\mu}_{1+a_{j_\rho}} \mbox{ \ \ for \ \ } \rho =1,...,r,$$
$$v_{2\rho -1} := \lambda _{1+j_{\rho -1}} - s, \ \ v_{2\rho} := \lambda_{j_\rho} - s \mbox{ \ \ for \ \ } \rho \neq \frac{r}{2},$$
$$v_{r-1} := \lambda _{1+ j_{\frac{r}{2}-1}} - s, \ \ v_r := \lambda _{j_{\frac{r}{2}}-1} -s.$$
\bigskip

{\bf Lemma 5.3:} {\it The weights $u$ and $v$ belong to $X_0^+(2r)$ and moreover the defect $d(u)$ vanishes.}
\bigskip 

{\bf Proof.} The weight $u$ is the same like in the non-exceptional case in Lemma 5.1. The weight $v$ is almost the same up to $v_r$ where purity follows by
$$v_r + v_{r+1} = \lambda _{j_\frac{r}{2}-1} - s + \lambda _{1+j_\frac{r}{2}} - s = \lambda _\frac{m-1}{2} + \lambda _\frac{m+3}{2} -2s = w'+n-m-1-2s,$$
and by Proposition 3.1 since by Remark 2.4 we have $j_\frac{r}{2} = \frac{m+1}{2}$. As to the defect we show that $1 + a_{j_\frac{r}{2}} = a_{j_{1+\frac{r}{2}}}$ which implies
$u_r=u_{r+1}$. By Remark 2.4 the index $j=\frac{m-1}{2}$ is not a jump index for even $r$ , so by Lemma 2.5 we have
$$a_\frac{m-1}{2} = a_\frac{m+1}{2} (= a_{j_\frac{r}{2}}) = \frac{n-1}{2}.$$
Since we know that $a_{j_{1+\kappa}}=a_{1+j_\kappa}$ for any jump index $j_\kappa$ we get by $(2.3)$ for $j=\frac{m-1}{2}$ the identity
$$a_{j_{1+\frac{r}{2}}} = a_{1+j_\frac{r}{2}} = a_\frac{m+3}{2} = n - a_\frac{m-1}{2} = \frac{n+1}{2} = 1 + a_{j_\frac{r}{2}}$$
as required.

Again we can apply Proposition 4.1 to $\tilde{\mu}:=\check{u}$ for $u\in X_{sp}(2r)$ and $\tilde{\lambda}:= v+s \in X_0^+(2r)$, hence we get that $(u,v)$ satisfies $(4.2)$
if and only if we have the embeddings
$$\theta _i'(\tilde{\lambda}) - s \ \ \prec \ \ \theta _i(\check{\tilde{\mu}}) \mbox{ \ \ for \ \ } i=1,2 \leqno(5.3)$$
just like in all other previous cases.

In order to take care of the parity condition $(PC)_\epsilon$ in Proposition 2.1 we introduce the notation 
$$Z_\epsilon := (2 \Bbb N - \epsilon ) \cup -(2 \Bbb N - 1 - \epsilon ) \mbox{\ for \ } \epsilon = 0,1$$
and we put 
$$\Emb(\alpha ,\beta )_\epsilon := \{s \in \Emb(\alpha ,\beta ); s-\frac{1}{2}(n-m+w+w')+1 \in Z_\epsilon \}$$
and $$(M_\alpha \otimes M_\beta )_\epsilon ^{SL_r} := \bigoplus _{s \in \Emb(\alpha , \check{\beta})_\epsilon} T_r(s).$$
The preceding discussion now eventually emerges into
\bigskip

{\bf Theorem B:} {\it In the exceptional case the mapping $t \mapsto t+\frac{n-m}{2} - 1$ sets up a bijection
$$\mathrm{Crit}(\pi _\infty ,\sigma _\infty ) \longrightarrow \bigcap _{i=1,2} \Emb(\theta _i'(\tilde{\lambda}),\theta _i(\check{\tilde{\mu}}))_\epsilon  .$$}
\bigskip

The analogue of the corollary of Theorem A may be formulated as follows. 
\bigskip  

{\bf Corollary 5.2:} {\it In terms of the Tate modules  $T_r(s)$ we have 
$$\bigcap _{i=1,2}(M_{\theta _i'(\tilde{\lambda})} \otimes M_{\theta _i(\tilde{\mu})})^{SL_r}_\epsilon \cong \bigoplus _{t\in Crit} T_r(t+\frac{n-m}{2}-1).$$
}
\bigskip

We finish by comparing our results with previously studied cases in the literature.
\bigskip

{\bf Examples:}

{\it Case n=2,\ m=1:} For a classical newform $f$ of weight $k$ lifting to an adelic function on $GL_2(\Bbb A)$ with central character $\omega$ with infinity part $\omega _\infty (x) = x^k$ for $x>0$  will provide us
with an attached cohomological cuspidal automorphic representation $\pi \in \Coh(\mu )$ with $\pi _\infty = J(-k,k-1)\otimes sign ^\delta$ and $\mu = (k-1,1)$. For trivial $\sigma$ 
we have $\nu =0, a=1, r=1$ and $d=0$, hence $\tilde{\mu}=\mu$ and $\tilde{\lambda}=(0,0)$. So Theorem A yields the bijection
$$\mathrm{Crit} \longrightarrow \Emb(0 ,\check{\mu}), \ t \mapsto t - \frac{1}{2}$$ 
since $\theta _i'(\tilde{\lambda}) = \nu$ and $\theta _i(\tilde{\mu}) = \mu $ which is a general fact
in the case $m=n-1$ (see below). Since for the completed $L$-function of the newform we have
$$L(\pi, s) = L(f,s + k - \frac{1}{2}),$$  
the  critical numbers are those studied  by  Shimura [10]. 

\bigskip  

{\it Case n=m=2:} For two classical newforms $f$ and $g$ of respective weight $k$ and $l$ with $k>l$ we may as in the previous case attach cohomological cuspidal automorphic representations $\pi$ and $\sigma$.
Since $a=(1,1)$ and $r=1$ 
we  get $\lambda =(l-1,0)$ and  $d = \lambda _1 - \lambda _2 = l-1,$ hence we have $\tilde{\mu} =(k-l,1), \tilde{\lambda} = (l-1,l-1), \theta _i'(\tilde{\lambda})=l-1$ 
and $\theta _i(\tilde{\mu}) = \tilde{\mu} = (k-l,1)$. We find
$$M_{\theta _i'(\tilde{\lambda})} \otimes M_{\theta _i(\tilde{\mu})} \cong M_{(k-1,l)} \cong \bigoplus _{s=l}^{k-1} T(s) \cong \bigoplus _{t\in \mathrm{\mathrm{Crit}}(\pi _\infty ,\sigma _\infty )} T(t-1).$$
Again we encounter the critical numbers studied by Shimura [loc.cit.] for  the associated L-function attached to the two forms $f$ and $g$. 

\bigskip

{\it Case m=n-1:}  We revisit the discussion in [7] where we assumed that $\\Emb(\nu ,\check{\mu})$ is non-empty. In that situation we have $a_j=j$, $j_\kappa = \kappa$, $r=m$ and $\lambda = \nu$. Moreover we get $d=0$, 
$\theta _i'(\tilde{\lambda }) = \nu$ and $\theta _i(\tilde{\mu }) = \mu$ so by Theorem A we have the bijection
$$\mathrm{\mathrm{Crit}}(\pi _\infty ,\sigma _ \infty ) \longrightarrow \\Emb(\nu ,\check{\mu}), \ \ t \mapsto t - \frac{1}{2}$$
and
$$(M_\mu \otimes M_\nu )^{SL_{n-1}} \cong \bigoplus _{t \in \mathrm{Crit}} T_{n-1}(t - \frac{1}{2}).$$

\bigskip

{\it Case $n=3,\ m=1$ The Jacquet-Gelbart-Lift:} Each cuspidal automorphic representation $\pi$ of $GL_2$ attached to a newform $f$ of weight $k$ as in the first example
gives rise to an automorphic representation $\Pi$ of $GL_3$ by the work of Gelbart and Jacquet [3]. The infinity component of $\Pi$ is 
$$ \Pi _\infty \cong J(0;(2(k-1),0,-2(k-1))) \otimes sgn$$
(see for instance [9] (1.8)). Hence for $\mu := (k-2,0,2-k)$ we find that $\Pi$ belongs to $\Coh(\mu )$ by [1] proof of Lemme 3.14. So we can apply Theorem B to the pair
$(\Pi _\infty , \sigma _\infty )$ for the trivial $GL_1$-representation $\sigma = {\bf 1}$. Since $\lambda = \nu =0$ and $d=0$ we have $\tilde{\mu} = (k-1,2-k)$ and $\tilde{\lambda} =(0,0)$. 
So we get with $\theta _i'(\tilde{\lambda})=0$ and $\theta _i(\tilde{\mu})= \tilde{\mu} = (k-1,2-k)$
$$\mathrm{Crit}(\Pi _\infty ) = \{2-k,...,k-1\}\cap (\{1,3,5,...\}\cup \{0,-2,-4,...\})$$
as considered in Lemma 2.1 of [9]. Moreover for any größencharacter $\xi$ with infinity component $\xi _\infty = {\mathrm sgn}\otimes |\cdot |_{\Bbb R}$ we have a cohomological $GL_1$-representation
$\sigma \in \Coh (\nu )$ for $\nu = -1$ hence with the attached $\lambda = -1$, $\tilde{\lambda} = (-1,-1)$, $\theta _i'(\tilde{\lambda}) = -1$, $\theta _i(\tilde{\mu}) = \tilde{\mu} = (k-1,2-k)$   again by Theorem B we find
$$\mathrm{Crit}(\Pi _\infty ,\sigma _\infty ) = \{1-k,...,k-2\}\cap (\{1,3,5,...\}\cup\{-2,-4,...\})$$ 
in accordance with Lemma 2.1 of [9] since we have
$$L(\Pi _\infty \otimes {\mathrm sgn}, s+1) = L(\Pi _\infty , \sigma _\infty , s).$$ 

\bigskip

{\bf References}

[1] {\it Clozel, L.:} Motifs et formes automorphes: Applications du principe de fonctorialité, in {\it Automorphic Forms, Shimura Varieties and L-Functions,} Perspectives in Mathematics, Vol.10
(Academic Press 1990) 77-159.

[2] {\it Deligne, P.:} Valeurs de fonctions $L$ et périodes d'intégrales. Proc. Symp. Pure Math. {\bf 33}, II (1979), 313-342. 

[3] {\it Gelbart, S., Jacquet, H.:} A relation between automorphic representations of $GL(2)$ and $GL(3)$. Ann. Sci. Ec. Norm. Super.,IV. Ser. {\bf 11} (1978), 471-542.             

[4] {\it Jacquet, H., Piatetski-Shapiro, I. I., Shalika, J. A.:} Rankin-Selberg convolutions. Am. J. Math. {\bf 105} (1983), 367-464. 

[5] {\it Januszewski, F.:} On Deligne's conjecture on special values of $L$-functions. Habilitationsschrift Karlsruhe KIT (2015)

[6] {\it Knapp, A.:} Local Langlands correspondence, the archimedean case. in Proc. Symp. Pure Mathematics, {\bf 55}(2), AMS 1994.

[7] {\it Kasten, H., Schmidt, C. G.:} On critical values of Rankin-Selberg convolutions. Int. J. Number Theory {\bf 09} (2013), 205-256.

[8] {\it Mahnkopf, J.:} Cohomology of arithmetic groups, parabolic subgroups and the special values of L-functions on $GL_n$. J. Inst. Math. Jussieu {\bf 4} (2005), 553-637.

[9] {\it Schmidt, C. G.:} P-adic measures attached to automorphic representations of $GL(3)$. Invent. math. {\bf 92}, (1988), 597-631. 

[10] {\it Shimura, G.:} On the periods of modular forms. Math. Annalen {\bf 229} (1977), 211-221.
\bigskip

Claus Günther Schmidt, Karlsruher Institut für Technologie, Institut für Algebra und Geometrie, Kaiserstraße 89-93,
76133 Karlsruhe, Germany, claus.schmidt@kit.edu

\end{document}